%% file: PSCCfullpaper.tex
\documentclass{IEEEtran4PSCC}
\ifCLASSINFOpdf
   \usepackage[pdftex]{graphicx}
\else
   \usepackage[dvips]{graphicx}
\fi
\usepackage[cmex10]{amsmath}
\makeatletter
\let\old@ps@headings\ps@headings
\let\old@ps@IEEEtitlepagestyle\ps@IEEEtitlepagestyle
\def\psccfooter#1{%
    \def\ps@headings{%
        \old@ps@headings%
        \def\@oddfoot{\strut\hfill#1\hfill\strut}%
        \def\@evenfoot{\strut\hfill#1\hfill\strut}%
    }%
    \def\ps@IEEEtitlepagestyle{%
        \old@ps@IEEEtitlepagestyle%
        \def\@oddfoot{\strut\hfill#1\hfill\strut}%
        \def\@evenfoot{\strut\hfill#1\hfill\strut}%
    }%
    \ps@headings%
}
\makeatother

\psccfooter{%
        \parbox{\textwidth}{\hrulefill \\ \small{23rd Power Systems Computation Conference} \hfill \begin{minipage}{0.2\textwidth}\centering \vspace*{4pt} \includegraphics[scale=0.06]{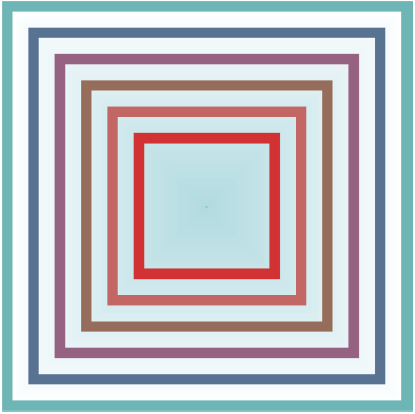}\\\small{PSCC 2024} \end{minipage} \hfill \small{Paris, France --- June 4 -- 7, 2024}}%
}

\usepackage{multirow}
\usepackage{mathtools}

\usepackage{amssymb}
\usepackage{rotating,stackengine,scalerel}
\usepackage{siunitx,booktabs,array}

\newcommand\wye{\scalerel*{\stackengine{-1pt}{%
  \rotatebox[origin=c]{30}{\rule{10pt}{.9pt}}\kern-1pt%
  \rotatebox[origin=c]{-30}{\rule{10pt}{1.3pt}}}{%
  \rule{.9pt}{10pt}}{O}{c}{F}{F}{S}}{\Delta}}

\begin{document}
%

\title{Improved Algebraic Inverter Modelling for Four-Wire Power Flow Optimization}

\author{
\IEEEauthorblockN{Rahmat Heidari}
\IEEEauthorblockA{Energy Business Unit, CSIRO \\
Newcastle, Australia}
\and
\IEEEauthorblockN{Frederik Geth}
\IEEEauthorblockA{GridQube \\
Brisbane, Australia}
}
\input{newcommands.tex}


\maketitle


\begin{abstract}
This paper discusses the modeling of inverters used in distributed energy resources in steady state.
Modeling the interaction between distribution grids and inverter-based resources is crucial to understand the consequences for the network's operational and planning processes. 
This work highlights the limitations of existing models and emphasizes the need for better representations of inverters and their control laws in decision-making contexts. 
Improved steady-state grid-following and grid-forming inverter models are presented, including both three-leg and four-leg converter variants. 
The advantages of these improved models in mathematical optimization contexts are showcased by investigating the power quality improvement capabilities of the inverters.
Numerical studies integrating the proposed inverter models in a four-wire unbalanced optimal power flow engine are presented, and trade-offs between modeling detail and computational intensity are illustrated.
\end{abstract}

\section{Introduction}
Driven by the energy transition, distributed energy resources (DERs) are common-place today in distribution networks.
DERs use power electronic converters to generate  voltage and/or current in synchronization with the network. 
Modelling of converters in simulation and optimization engines supports power engineers in addressing operational and planning challenges that stem from high numbers of DERs.

The most detailed existing models of inverters are in the electromagnetic transient (EMT) domain \cite{Lara_2024, Ghazavi9678003}.  
EMT models provide a thorough description of the inverter's internals circuitry including representations for power electronic switches,  topology, and control and protection systems.
Such models however are computationally expensive.
Reference frame transformations \cite{Rourke8836094} and average models \cite{Chiniforoosh5440017} via dynamic phasors are then used to derive equivalent or simplified formulations for which larger time steps could be used. Quasi-static phasor formulations that algebraically represent system dynamics as discrete changes between steady-state operating points are
used to study low-frequency phenomena ranging from inertial response to frequency regulation \cite{Demiray_thesis}.
Steady state studies such as load flow and short-circuit studies \cite{Cunha2021, 9916672} are less computationally demanding and enable engineers to perform extensive studies with much longer time horizons. 

\subsection{Steady State Inverter Models and Controllability}
Three-phase inverters exist with various circuit topologies. 
For example, they can have either three or four legs (with a leg being a pair of power electronic switches), where the optional fourth leg would connect to the neutral. 
Such topology differences impact the ability to control current in the phases independently, which in turn limits the feasible set points in terms of active and reactive power.


These kinds of details in inverter designs also impact the ability to deliver power quality services such as phase unbalance compensation \cite{Tomislav14010117}. 
Therefore, they should also be represented in decision-making contexts where phase unbalance is a  concern.
For instance, the quantification of hosting capacity or dynamic operating envelopes will be different depending on inverter control capabilities. 

Furthermore, the inverter \emph{control laws} (a.k.a. control modes) limit or restrict the degrees of control freedom \cite{GEBBRAN2021116546}.
For instance, Volt-var/Watt control modes in PV inverters have become mandatory for new installations in places such as Australia.
In the context of microgrids without synchronous generators, three-phase inverters should incorporate control laws to push the three-phase voltage phase angles apart to about 120 degrees, thereby (partially) emulating the behavior of synchronous generators. 
Detailed modelling of these control laws in steady state improves power system studies by providing accurate information about feasible regions as well as delivering regulated outputs that can be dispatched.
Modelling grid forming inverters is also important in the context of microgrid (or weak grid) dispatch or design optimization questions.

Devices such as solar or battery inverters have built-in control \emph{flexibility}. 
The maximum power point tracking (MPPT) functionality can be deviated from (i.e. curtailment), and even \emph{needs} to be deviated from in real-world situations.
For instance, it becomes common to oversize the panels (kWpeak) significantly relative to the rating of the inverter, as this can boost generation volumes in the morning and evening when power is often more valuable on the market.
On the sunniest days of the year though, this means that the MPPT may be restricted by the power rating of the inverter itself, thereby necessarily curtailing a bit of the capturable solar power. 
Similarly, reactive power control (mostly free except for a slight increase in inverter losses), can be exploited to support or suppress grid voltage magnitude.
These control laws have been developed as distributed control and/or fallback mechanisms, to work in conjunction with centralized control paradigms. As they are meant to be active continuously, we need to take the response of these mechanisms into account when developing centralized control algorithms such as unbalanced optimal power flow.
In distribution systems, deviation from fully balanced voltage and current phasors implies a higher ratio of losses to power being transferred, i.e. avoidable losses, which reduces network capacity utilization. 

\subsection{Scope and Contributions}
In this work, we exclude diode-based rectifiers or rectifier-only approaches, i.e. focus on power electronic converters with active switching. Furthermore, we assume a fixed electrical frequency and leave variable frequency and harmonics for future work. We also focus on 3- and 4-leg inverters and limit the discussion of single-phase inverters as a mathematical edge case of a 4-leg one.

Noting the previously discussed limitations in the literature, the aim of this paper is to:
\begin{itemize}
    \item perform a review  of algebraic inverter models and formulate improvements on the existing ones;
    \item derive (algebraic) nonlinear programming models for both the inverters and possible control laws, and perform numerical experiments through \textsc{Ipopt} \cite{Wachter2006};
    \item demonstrate the integration of these models into  unbalanced optimal power flow engines, up to four-wire \cite{CLAEYS2022108522}; 
    \item illustrate how inverter model approximations can lead to sub-par decisions.
\end{itemize}

The paper is structured as follows.
Section \ref{sec_lit} reviews the state of the art on steady-state inverter modeling and control.
Next, section \ref{sub_inverter_four_wire} proposes a constructive mathematical framework that enables the modeling of a variety of inverter topologies and control goals, compatible with unbalanced optimal power flow models.
Moreover, section \ref{sec_studies} showcases numerical optimization studies for a variety of inverter models stemming from the previously proposed framework.
Finally, section \ref{sec_conclusions} summarizes results and concludes the paper.

\section{Review of Inverter Models and Controls} \label{sec_lit}
We separate the parts of the models that relate to topology (P/Q space characterization) from those of control laws (volt-var, unbalance compensation etc).


\subsection{Inverters in Unbalanced Optimal Power Flow}
Accurate representation of inverter control strategies improves operation and optimization of distribution systems and DERs such as battery storage systems \cite{geth2015, Geth2020, aaslid2020}.
Battery storage systems delivering for network support services are considered an opportunity to manage congestion in distribution networks, which have become commonplace in different countries \cite{Beckstedde2023}. 
%
Specifically, network services such as voltage regulation benefit from more detailed models of the network physics, such as including mutual inductance between conductors carrying unbalanced currents, as well as inverter abilities to control current independently in the phases. 
Active filters are inverters with specialized control loops that can support inter-phase power exchange, unbalanced reactive power set points, and harmonic compensation.
While phase unbalance can be compensated by active filters \cite{Singh1999}, i.e. inverters without a prime mover (no solar/batteries), inverters with a prime mover provide even more flexibility that can be harvested for DER orchestration.

The ability of inverters to control voltage magnitude depends on the resistance and reactance of the network. 
Debrabandere et al. in \cite{Brabandere2007} provide a series of intuitive visualizations detailing controllability for networks with different ratios of resistance to reactance.

\begin{figure}[tbh]
  \centering
    \includegraphics[width=0.4\columnwidth]{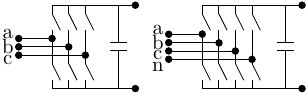}
  \caption{3-Leg vs. 4-Leg inverter power electronic circuit topologies.}  \label{fig_differentlegs}
\end{figure}

Inverter topology plays an important role in phase unbalance control in the network.  Stuyts et al. \cite{Stuyts2018} study the unbalance compensation capability of 3-leg and 4-leg inverters (Fig. \ref{fig_differentlegs}), connected in parallel with 1-phase, 2-phase and 3-phase loads. 
Both topologies are capable of compensating negative sequence component of the unbalanced current. However, the main difference between 3- and 4-leg converters is that the 4-leg converters can compensate zero sequence current which in the 3-leg case must be supplied by the external system, i.e. the grid.
This provides a distinctive advantage for the 4-leg topology.
A drawback for zero-sequence current compensation in 4-leg converters is that often it requires oversizing of the neutral conductor in the converter.

\subsection{Inverter Control Degrees of Freedom}

Grid following inverters typically use a phase-locked-loop and a current control loop to achieve fast control of the inverter’s output currents. 
We do not capture models for these control loops in the steady state model developed in this paper, as they are effects which play out on a time scale of less than 20\,ms\footnote{at 50\,Hz, or the equivalent period at 60\,Hz.}. 
Grid-following inverters typically regulate a constant power output. 
However, they rely on an external voltage source to provide the voltage and frequency references. 
Control laws applied to the grid following inverters include  Volt-var and Volt-Watt control modes \cite{Lusis2020}, constant power factor, constant reactive power, power sharing, and power rate limits, which the latter is a dynamic control and out of scope for this paper.

Grid-forming inverters, on the other hand,  require a primary source of power (a.k.a. prime mover) and can autonomously establish voltage magnitude and frequency for a network. 
Control systems for grid-forming inverters often aim to emulate the behavior of synchronous generators which produce a balanced internal voltage phasor of which the magnitude cannot be changed. Grid forming controls include emulating a positive sequence voltage source, droop controls such as voltage-frequency droops, and synchronous machine behavior replication.
For instance  inertia emulation \cite{Rathnayake9513281}  is accomplished through the swing equation
\begin{equation}
    \label{eq:swing}
    \frac{2H}{\omega_s} \frac{d^2 \theta}{dt^2} = P_m - P_e
\end{equation}
where $H$ is the inertia constant, $\omega_s$ is the steady-state frequency, $\theta$ is the rotor angle with respect to the synchronously rotating reference frame, $t$ denotes time, and $P_m$ and $P_e$ are mechanical and electrical powers respectively.
As the equation \eqref{eq:swing} is a differential equation, in keeping aligned with the steady-state scope of this work, inertia emulation and frequency droop are \emph{not} included in the  models proposed in the upcoming sections, but the electrical angle $\theta$ is maintained as a variable.

Most of the literature on models for grid-forming inverters assumes balanced, inductive networks \cite{Rathnayake9513281}.
In this work we focus on 1) distribution networks with 2) unbalanced operation and 3) and fast-acting droops.
Details of control strategies implemented in  commercially-available inverters are generally not well-described in spec sheets or brochures, so we have to rely on reports that detail observed limitations of such solutions.
The \emph{Research roadmap on grid-forming inverters} \cite{osti_1721727} states:
``no commercially available three-phase grid-forming inverters independently balance phases by providing zero- and negative-sequence currents (both as a result of the current $dq$ control paradigm as well as the tendency of commercial inverters to come in three-leg, instead of four-leg, configurations as a cost-saving measure)''.
Furthermore, ``traditional three-phase grid-following controls will not provide zero- or negative-sequence currents  \cite{osti_1721727}''.

The potential of smarter controls is also recognized.
Merritt et al. in \cite{Merritt2017} state that grid-forming inverters can reduce the negative-sequence voltage to a negligible level by controlling the positive- and negative-sequence voltages separately.
Kim et al. \cite{Kim2022} state that although inverters may be equipped with a negative-sequence voltage controller, they often cannot supply a large amount of negative-sequence current due to their low rated current compared to synchronous machines.

Cunha et al. \cite{Cunha2021} derive positive- and negative-sequence equivalent circuits of three-phase three-leg grid-forming inverters for steady state analysis in $dq$- and $\alpha\beta$-frames. 
The steady state models implemented in  \textsc{OpenDSS} are then validated against time-domain models in PSCAD. 
An method to model the current limiter is proposed, but the authors address that representing the current limiter in $abc$-frame has better performance considering unbalanced scenarios such as during unbalanced loads and faults, which is aligned with the modelling framework considered in this work.

The inverter model in GridLab-D is detailed in \cite{Wei2021} and that of  \textsc{OpenDSS} in \cite{opendss_inverter_controls}.
PowerModelsDistribution also has an inverter model for battery storage built in \cite{FOBES2020106664}.
Aaslid et al. develop an empirical model for inverters \cite{aaslid2020}, that is compatible with nonlinear programming frameworks. 
Savasci et al. propose mathematical models to add Volt-var/Watt control laws for inverters to unbalanced OPF problems \cite{Adedoyin9796576}, but end up linearizing the network physics, as their approach relies on the encoding of the piecewise functions using binary variables.
Quiertant et al. \cite{Quiertant10407677} study active voltage management using LV STATCOM or PV inverter  Volt-var control using nonlinear programming techniques.

In this work we aim to take a step towards addressing ``the need to begin developing appropriate models for existing simulation tools as well as enhanced modeling and simulation tools'', as identified in \cite{osti_1721727}.
We observe that there is a gap in the literature on the modeling of inverters when it comes to the representation of the degrees of freedom in the control systems.
Modeling of Volt-var/Watt control loops in optimization algorithms is the one that has received a bit of attention, however papers have rarely discussed how to establish models with inter-phase power exchange and the nuances between 3- and 4-leg designs. 

\section{Steady-State Inverter Models up to  Four-Leg} \label{sub_inverter_four_wire}
\begin{figure}[tbh]
  \centering
    \includegraphics[width=1\columnwidth]{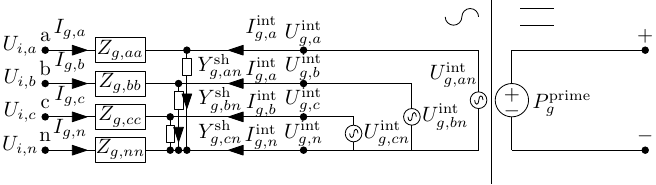}
  \caption{Equivalent circuit of the four-leg inverter model with mathematical symbols, and current directions indicated. }  \label{fig_inverter}
\end{figure}

Fig. \ref{fig_inverter} defines the most important scalar voltage and current variables used in the inverter model, both on the bus it is connected to, and the internal voltages and the prime mover. 

We represent up to four nodes per bus, $\mathcal{N} = \{a,b,c,n \}$, which are composed of the phases $\mathcal{P} \in \{a,b,c \}$ and the neutral $n$. 
We use a bold typeface to indicate a vector over the nodes and/or phases, e.g. for complex voltage $\mathbf{U}$, current $\mathbf{I}$ and  power $\mathbf{S}$. 
Sometimes, it is practical to refer to an extract of a vector only defined over the phases, from a vector that also includes the neutral. 
For instance, the phase current vector derived from the length-four node vector $\mathbf{I} := \mathbf{I}[\mathcal{N}]$  is written as $\mathbf{I}[\mathcal{P}]$.

In the next sections we first develop the model for a 4-leg grid following inverter, and then illustrate the further adaptation to a grid-forming inverter. 
Similarly, we derive the model for the 3-leg inverter from the 4-leg one, as a restriction.

\subsection{Grid-Following  Inverter (GFL)}
\label{sec:GFL}
We first provide the more general four-leg representation of inverters before introducing detailed model of grid-following and grid-forming inverters.

\subsubsection{Voltage}
At the point of common coupling (PCC), bus phase-to-ground and phase-to-neutral (`$\wye$') voltages are,
\begin{IEEEeqnarray}{C}
    \label{eq:Ui}
    \mathbf{U}_i = 
    \begin{bmatrix} 
        U_{i,a} \\
        U_{i,b} \\
        U_{i,c} \\
        U_{i,n} \\
    \end{bmatrix},
    \label{eq:Ui_pn}
    \mathbf{U}^{\wye}_i = 
    \begin{bmatrix} 
        U_{i,an} \\
        U_{i,bn} \\
        U_{i,cn} \\
    \end{bmatrix}
    =
        \begin{bmatrix} 
        U_{i,a} - U_{i,n}  \\
        U_{i,b} - U_{i,n} \\
        U_{i,c}- U_{i,n} \\
    \end{bmatrix}.
\end{IEEEeqnarray}
The phase-to-phase voltages (`$\Delta$') are,
\begin{IEEEeqnarray}{C}
    \mathbf{U}^{\Delta}_i = 
    \begin{bmatrix} 
        U_{i,ab} \\
        U_{i,bc} \\
        U_{i,ca} \\
    \end{bmatrix}
    =
        \begin{bmatrix} 
        U_{i,a} - U_{i,b}  \\
        U_{i,b} - U_{i,c} \\
        U_{i,c}- U_{i,a} \\
    \end{bmatrix}. \label{eq_v_delta}
\end{IEEEeqnarray}

The inverter $g$ \emph{internal} phase-to-ground and phase-to-neutral voltages are,
\begin{IEEEeqnarray}{C}
    \label{eq:Ugint}
    \Ugint = 
    \begin{bmatrix} 
        U^{\text{int}}_{g,a} \\
        U^{\text{int}}_{g,b} \\
        U^{\text{int}}_{g,c} \\
        U^{\text{int}}_{g,n} \\
    \end{bmatrix},
    \label{eq:Ugint_pn}
    \mathbf{U}^{\text{\wye,int}}_{g} = 
    \begin{bmatrix} 
        U^{\text{int}}_{g,an} \\
        U^{\text{int}}_{g,bn} \\
        U^{\text{int}}_{g,cn} \\
    \end{bmatrix} =
    \begin{bmatrix} 
        U^{\text{mag}}_{g,an} \angle U^{\text{ang}}_{g,an} \\
        U^{\text{mag}}_{g,bn} \angle U^{\text{ang}}_{g,bn} \\
        U^{\text{mag}}_{g,cn} \angle U^{\text{ang}}_{g,cn} \\
    \end{bmatrix}.
\end{IEEEeqnarray}
%
%

\subsubsection{Current}
Conductor currents at the bus are,
\begin{IEEEeqnarray}{C}
    \label{eq:Ii}
    \mathbf{I}_g :=     \mathbf{I}_g[\mathcal{N}] = 
    \begin{bmatrix} 
        I_{g,a} \\
        I_{g,b} \\
        I_{g,c} \\
        I_{g,n} \\
    \end{bmatrix},
    \label{eq:Ii_p}
    \mathbf{I}_g[\mathcal{P}] = 
    \begin{bmatrix} 
        I_{g,a} \\
        I_{g,b} \\
        I_{g,c} \\
    \end{bmatrix}.
\end{IEEEeqnarray}
The inverter's \emph{internal} conductor currents are,
\begin{IEEEeqnarray}{C}
    \label{eq:Ig}
    \mathbf{I}^{\text{int}}_{g} = 
    \begin{bmatrix} 
        I^{\text{int}}_{g,a} \\
        I^{\text{int}}_{g,b} \\
        I^{\text{int}}_{g,c} \\
        I^{\text{int}}_{g,n} \\
    \end{bmatrix},
    \label{eq:Ig_p}
        \mathbf{I}^{\text{int}}_{g}[\mathcal{P}] = 
    \begin{bmatrix} 
        I^{\text{int}}_{g,a} \\
        I^{\text{int}}_{g,b} \\
        I^{\text{int}}_{g,c} \\
    \end{bmatrix},
\end{IEEEeqnarray}
which (in non-fault conditions) are expected to satisfy KCL,
\begin{IEEEeqnarray}{C}
    \label{eq:Ig_sum_0}
    \mathbf{1}^{\text{T}}  \mathbf{I}^{\text{int}}_{g} = 0. 
\end{IEEEeqnarray}
To protect the switches, e.g. IGBTs, against thermally-induced damage, the internal  currents are limited below $\mathbf{I}^{\text{rating}}_{g}$,
\begin{IEEEeqnarray}{C}
    \label{eq:Ig_limit}
    \mathbf{I}^{\text{int}}_{g}  \circ (\mathbf{I}^{\text{int}}_{g})^* \le \mathbf{I}^{\text{rating}}_{g} \circ \mathbf{I}^{\text{rating}}_{g}.
\end{IEEEeqnarray}
Note that this bound represents the current limiter in inverter control structure. 

We assume series-shunt output filter in series with the power-electronic circuit. 
The series part is represented by $Z_{g,aa}, Z_{g,bb}, Z_{g,cc}, Z_{g,nn}$. As this path is electrically short, we do not model \emph{mutual} impedance, just \emph{self} impedance.
The series filter part is typically implemented as inductances. 
The shunt part, typically implemented as high-voltage capacitors, is typically connected in wye configuration between the phases and optionally to the neutral wire. 
For the three-leg case, the star point is left floating.
The inverter shunt current is, 
\begin{IEEEeqnarray}{C}
    \label{eq:Ig_shunt}
    \mathbf{I}^{\text{sh}}_{g} = 
    \begin{bmatrix}
        I^{\text{sh}}_{g,a}  &
        I^{\text{sh}}_{g,b} &
        I^{\text{sh}}_{g,c}  &
        I^{\text{sh}}_{g,n}  
        \end{bmatrix}^{\text{T}}
=     \mathbf{Y}^{\text{sh}}_g \circ \mathbf{U}_i
\end{IEEEeqnarray}
The current vectors in the inverter model satisfy,
\begin{IEEEeqnarray}{C}
    \label{eq:Ig_Ish_Ii}
    \mathbf{I}_{g} + \mathbf{I}^{\text{int}}_{g} = 
    \mathbf{I}^{\text{sh}}_{g}.
\end{IEEEeqnarray}


\subsubsection{Power}
The external inverter power at the PCC is,
\begin{IEEEeqnarray}{C}
    \label{eq:Si}
    \mathbf{S}_{g} = 
        \mathbf{P}_{g} +  j \mathbf{Q}_{g}
    = 
    \mathbf{U}_i  \circ  ( \mathbf{I}_g )^{*} \nonumber
    \\
    =
    \begin{bmatrix}
        S_{g,aa}  \\
        S_{g,bb} \\
        S_{g,cc}  \\            
        S_{g,nn}  \\            
    \end{bmatrix}
    =
    \begin{bmatrix}
        P_{g,aa}  + jQ_{g,aa} \\
        P_{g,bb}  + jQ_{g,bb} \\
        P_{g,cc}  + jQ_{g,cc} \\
        P_{g,nn}  + jQ_{g,nn} \\   
    \end{bmatrix}.
\end{IEEEeqnarray}
The internal inverter power is,
\begin{IEEEeqnarray}{C}
    \label{eq:Sgint}
    \Sgint =  \Pgint +    j     \Qgint
    = 
    \Ugint \circ ( \mathbf{I}^{\text{int}}_{g} )^{*}.
\end{IEEEeqnarray}
The inverter filter complex power losses are,
\begin{IEEEeqnarray}{C}
    \label{eq:Sgint_Si}
 \mathbf{S}_{g}  +   \Sgint 
    = \mathbf{S}^{\text{loss}}_g.
\end{IEEEeqnarray}
Phase-to-neutral ($\wye$) external  power is,  
\begin{IEEEeqnarray}{C}
    \label{eq:Si_pn}
    \mathbf{S}^{\wye}_{g} = 
    \begin{bmatrix}
        S^{\wye}_{g,an} &
        S^{\wye}_{g,bn} &
        S^{\wye}_{g,cn} 
    \end{bmatrix}^{\text{T}} =
    \mathbf{U}^{\wye}_i  \circ ( \mathbf{I}_{g}[\mathcal{P}] )^{*} .
\end{IEEEeqnarray}

The complex power generated by the inverter's internal voltage source is summed across the phases,
\begin{IEEEeqnarray}{C}
    \label{eq:Sgint_sum}
     \mathbf{1}^{\text{T}}  \Sgint = {P}^{\text{int}}_{g} + j {Q}^{\text{int}}_{g}.
\end{IEEEeqnarray}
We assume that the inverter's internal active power is provided by a primary source, e.g. PV, wind, or battery storage,
\begin{IEEEeqnarray}{C}
    \label{eq:Pgint_pprimary}
    {P}^{\text{int}}_{g} = {P}^{\text{prime}}_{g},
\end{IEEEeqnarray}
with feasible states for prime mover control in a yet-to-be-specified set, i.e. ${P}^{\text{prime}}_{g} \in \mathcal{P}_g$ . These  sets may model effects such as standy losses of the inverter itself, deviation from MPPT, battery charging, etc. 
Note that we do not define voltage or current variables associated with the prime mover, and leave that for future work.
The reactive power balanced is controllable as well, 
\begin{IEEEeqnarray}{C}
    \label{eq:Qgint_slack}
    {Q}^{\text{int}}_{g} = {Q}^{\text{slack}}_{g},
\end{IEEEeqnarray}
with feasible states for reactive power control in a yet-to-be-specificed set, i.e. ${Q}^{\text{slack}}_{g} \in \mathcal{Q}_g$.


\subsubsection{Mathematical Reduction to a 3-Leg Model} \label{sec:GFL_3w}
To the previously-discussed system of constraints, we add the additional constraint:
\begin{IEEEeqnarray}{C}
    \label{eq:Ign_0}
 I_{g,n} =  I^{\text{int}}_{g,n} = 0.
\end{IEEEeqnarray}
And we set series/shunt impedance values for the neutral to be open circuits, so $U^{\text{int}}_{g,n}$ is floating and can be eliminated.
In combination with KCL \eqref{eq:Ig_sum_0} this  forces the phase currents to add to 0.
In implementations it is generally deemed prudent to eliminate the redundant variables, instead of assigning a fixed value. 
The expressions for internal phase-to-neutral voltage/power/current variables become meaningless as well. 

\subsection{Grid-Forming Inverter (GFM)}
We now cast Grid-Forming Inverters as a restriction of GFL, i.e. develop additional constraint sets.
A GFM inverter behaves as a power source which also controls the internal voltage source.
For instance, the internal voltage source being three-phase balanced,
\begin{IEEEeqnarray}{C}
    \label{eq:GFM_Umag}
    U^{\text{mag}}_{g,an} = U^{\text{mag}}_{g,bn} = U^{\text{mag}}_{g,cn}, \\
    \label{eq:GFM_Uangle1}
    U^{\text{ang}}_{g,an} - U^{\text{ang}}_{g,bn}= U^{\text{ang}}_{g,bn} -U^{\text{ang}}_{g,cn} = {2\pi}/{3} .
\end{IEEEeqnarray}
Note that the third permutation, $U^{\text{ang}}_{g,cn} - U^{\text{ang}}_{g,an} ={2\pi}/{3}$ is redundant w.r.t. \eqref{eq:GFM_Uangle1}. 
Also note that in a synchronous generator, the voltage magnitude of the internal voltage source is also fixed, however, this does not \emph{need} to be enforced in the inverter control system.
This combination of constraints can be equivalently cast as the internal voltage source being purely positive-sequence.

\subsection{Control Laws}
We now list a variety of control laws that have been discussed in the literature using the notation introduced previously. The control laws presented here are a subset of what is available in the existing literature \cite{Braun6733334}.

\subsubsection{Droop control}
Inverter voltage magnitude can vary with linear functions of active and reactive power (without breakpoints), 
\begin{IEEEeqnarray}{C}
    \notag
|\mathbf{U}^{\wye}_i|  - \mathbf{U}^{\text{setpoint}}_g  = D^{\text{q}} (\mathbf{Q}_{g}  - \mathbf{Q}^{\text{setpoint}}_{g})
    \\
    \label{eq:control_droop}
    \qquad \qquad \qquad
    + D^{\text{p}} ( \mathbf{P}_{g} - \mathbf{P}^{\text{setpoint}}_{g}),
\end{IEEEeqnarray}
where $\mathbf{U}^{\text{setpoint}}_g, \mathbf{Q}^{\text{setpoint}}_{g}, \mathbf{P}^{\text{setpoint}}_{g}$ is are the set points for voltage, active and reactive power,  $|\mathbf{U}^{\wye}_i|$ is the vector of phase-to-ground voltage magnitudes, and $D^p$ and $D^q$ are droop control coefficients.

\subsubsection{Volt-var and Volt-Watt control}
Volt-var $f^{VV}(\cdot)$ and Volt-Watt $f^{VW}(\cdot)$ functions today are common-place functionalities offered by mass-market solar inverters. 
We note that these nonlinear functions, typically defined by a set of break points of a piece-wise linear function, may differ between jurisdictions. 
We define variants through equations \eqref{eq:control_4wire_Pan}-\eqref{eq:control_ODSS_Qag} in Table~\ref{tab:vv_vw} in terms of phase-to-neutral, phase-to-phase and phase-to-ground voltages, either per phase, or across the phases.
The averaged phase-to-ground variants  allow us validate w.r.t. the inverter control functionality offered in \textsc{OpenDSS}.

\begin{table}
    \caption{Volt-var and Volt-Watt Control Variations.}
    \label{tab:vv_vw}
    \setlength{\tabcolsep}{2pt}
    \def\arraystretch{1.2}
    \begin{tabular*}{0.48\textwidth}{
        @{\extracolsep{\fill}}
        l
      >{$\displaystyle}c<{\vphantom{\sum_{1}{N}}$}
      >{\refstepcounter{equation}(\theequation)}r
      @{}
    }
        \toprule
        Variation & \multicolumn{1}{l}{Control} & \multicolumn{1}{r}{}
        \\
        \hline
        \multirow{2}{4.4em}{Phase-to-neutral} &
        \multicolumn{1}{l}{$P_{g,an} =  f^{VW}(U^{\text{mag}}_{g,an})$} & \label{eq:control_4wire_Pan}
        \\
         & \multicolumn{1}{l}{$Q_{g,an} =  f^{VV}(U^{\text{mag}}_{g,an})$} & \label{eq:control_4wire_Qan}
         \\ \hline
        \multirow{2}{4.2em}{Phase-to-phase} & 
        \multicolumn{1}{l}{$P_{g,ab} =  f^{VW}(U^{\text{mag}}_{g,ab}/\sqrt{3} )$} & \label{eq:control_3wire_Pab}
        \\
         & \multicolumn{1}{l}{$Q_{g,ab} =  f^{VV}(U^{\text{mag}}_{g,ab}/\sqrt{3} )$} & \label{eq:control_3wire_Qab}
        \\ \hline
        \multirow{2}{4.2em}{Phase-to-phase averaged} & 
        \multicolumn{1}{l}{$P_{g,ab}=P_{g,bc}=P_{g,ca} =  f^{VW}\left(\frac{ U^{\text{mag}}_{g,ab} + U^{\text{mag}}_{g,bc} + U^{\text{mag}}_{g,ca}}{3\sqrt{3}} \right)$} & \label{eq:control_3wire_Pab_avg}
        \\
        & \multicolumn{1}{l}{$Q_{g,ab}=Q_{g,bc}=Q_{g,ca} =  f^{VV}\left(\frac{ U^{\text{mag}}_{g,ab} + U^{\text{mag}}_{g,bc} + U^{\text{mag}}_{g,ca}}{3\sqrt{3}} \right)$} & \label{eq:control_3wire_Qab_avg}
        \\\hline
         \multirow{2}{4.4em}{Phase-to-ground averaged} & 
        \multicolumn{1}{l}{$P_{g,a}= P_{g,b}= P_{g,c} =  f^{VW}\left(\frac{ U^{\text{mag}}_{g,a} + U^{\text{mag}}_{g,b} + U^{\text{mag}}_{g,c}}{3} \right)$} & \label{eq:control_3wire_Pag_avg}
         \\
         & \multicolumn{1}{l}{$Q_{g,a}= Q_{g,b}=Q_{g,c} =  f^{VV}\left(\frac{ U^{\text{mag}}_{g,a} + U^{\text{mag}}_{g,b} + U^{\text{mag}}_{g,c}}{3} \right)$} & \label{eq:control_3wire_Qag_avg}
         \\ \hline
         \multirow{2}{4.4em}{Phase-to-neutral averaged} & 
        \multicolumn{1}{l}{$P_{g,an}= P_{g,bn}= P_{g,cn} =  f^{VW}\left(\frac{ U^{\text{mag}}_{g,an} + U^{\text{mag}}_{g,bn} + U^{\text{mag}}_{g,cn}}{3} \right)$} & \label{eq:control_3wire_Pan_avg}
         \\
         & \multicolumn{1}{l}{$Q_{g,an}= Q_{g,bn}=Q_{g,cn} =  f^{VV}\left(\frac{ U^{\text{mag}}_{g,an} + U^{\text{mag}}_{g,bn} + U^{\text{mag}}_{g,cn}}{3} \right)$} & \label{eq:control_3wire_Qan_avg}
         \\ \hline
        \multirow{2}{4.4em}{Phase-to-ground} & 
        \multicolumn{1}{l}{$P_{g,a} =  f^{VW}(U^{\text{mag}}_{g,a} )$} & \label{eq:control_ODSS_Pag}
        \\
         & \multicolumn{1}{l}{$Q_{g,a} =  f^{VV}(U^{\text{mag}}_{g,a} )$} & \label{eq:control_ODSS_Qag}
         \\ 
         \bottomrule
    \end{tabular*}
\end{table}

\subsubsection{Power Sharing Control}
Four-leg inverters, can be restricted to balanced power output \cite{geth2015} on the phases,  
\begin{IEEEeqnarray}{c}
S^{\wye}_{g,an} = S^{\wye}_{g,bn} = S^{\wye}_{g,cn}.
\end{IEEEeqnarray}

\subsubsection{Constant Power Factor}
The inverter active power depends on the primary source power \eqref{eq:Pgint_pprimary} and therefore the inverter reactive power in constant power factor mode is controlled as
\begin{IEEEeqnarray}{C}
     \label{eq:control_const_pf}
    {Q}^{\text{int}}_{g} = {P}^{\text{int}}_{g} \tan(\arccos(\text{PF})),
\end{IEEEeqnarray}
where PF is the power factor.

\subsubsection{Sequence Current Control}

Inverters controls can be set up in symmetrical components instead of the phase coordinates. 
The transformation to symmetrical components is, 
\begin{IEEEeqnarray}{C}
    \label{eq:control_Ig_012_target}
    \mathbf{I}_{g}^{012} =
    \begin{bmatrix}
        I_{g,0} \\
        I_{g,1} \\
        I_{g,2} \\
    \end{bmatrix}
    =
  \frac{1}{3}
    \begin{bmatrix} 
        1 &  1  & 1 \\
        1 & \alpha   & \alpha^2 \\
        1 &  \alpha^2  & \alpha \\
    \end{bmatrix}
    \begin{bmatrix}
        I_{g,a} \\
        I_{g,b} \\
        I_{g,c} \\
    \end{bmatrix}
     = 
     \mathbf{T} \mathbf{I}_{g}[\mathcal{P}],
\end{IEEEeqnarray}
where $\alpha = e^{j 2\pi/3}$.
Note that this definition implies $ I_{g,0} =  -I_{g,n}/3$.
Limits on the magnitudes of the symmetrical component currents can be enforced,
\begin{IEEEeqnarray}{C}
    \label{eq:Ig_limit_012}
     \mathbf{I}_{g}^{012}   \circ ( \mathbf{I}_{g}^{012} )^* \le  \mathbf{I}_{g}^{012, \text{max}} \circ \mathbf{I}_{g}^{012, \text{max}}.
\end{IEEEeqnarray}
Pure positive sequence current control can be achieved by setting the upper bounds for negative and zero-sequence magnitude currents (i.e. the corresponding entries $\mathbf{I}_{g}^{012, \text{max}}$) to zero. 



\subsection{Phase Unbalance Control}


For a 4-leg inverter, the sequence voltages are conventionally defined w.r.t the phase-to-neutral voltages, 
\begin{IEEEeqnarray}{C}
    \label{eq_voltage_sequence}
 \mathbf{U}_{i}^{012} 
 =
 \begin{bmatrix}
        U_{i,0} \\
        U_{i,1} \\
        U_{i,2} \\
    \end{bmatrix}
    =
  \frac{1}{3}
    \begin{bmatrix} 
        1 &  1  & 1 \\
        1 & \alpha   & \alpha^2 \\
        1 &  \alpha^2  & \alpha \\
    \end{bmatrix}
    \begin{bmatrix} 
        U_{i,an} \\
        U_{i,bn} \\
        U_{i,cn} \\
    \end{bmatrix}
    =
     \mathbf{T}\mathbf{U}^{\wye}_i.
\end{IEEEeqnarray}
In networks without neutral, transformation $\mathbf{T}$ is applied to $  \mathbf{U}^{\Delta}_i$ instead. 

From the technical standards of the electricity distribution code review \cite{Gridcode2019}, the positive sequence voltage of an unbalanced  phasor cannot exceed $[-10, 10]\%$ of the nominal voltage $U^{\text{nom}}$,
\begin{align}
    \label{eq:Ui1_limits}
    & |U_{i,1}| \in [0.9  \cdot U^{\text{nom}}, 1.1  \cdot U^{\text{nom}}].
\end{align}
The negative sequence component is also constrained to $2\%$, that is,
\begin{align}
    \label{eq:Ui2_limit}
    & |U_{i,2}| \le 2 \% \cdot U^{\text{nom}}.
\end{align}
In unbalanced distribution systems, an interesting control objective is to minimize the unbalance. 
There are several definitions of unbalance applied to both voltage and current variables \cite{Bollen2002}. 
From \cite{NEMA1993} the  definition of the \emph{voltage unbalance factor} (VUF) is,
\begin{IEEEeqnarray}{C}
    \label{eq:VUF}
    \text{VUF} = |U_{i,2}| / |U_{i,1}| \le 2 \% .
\end{IEEEeqnarray}
More definitions of voltage unbalance are proposed in \cite{IEEE8291810}, considering phase and line voltages. 

\subsection{Feasible Sets for the Various Inverter Models}
\label{sec_feasibility_regions}

Tables~\ref{table:inverer_constraints} and \ref{table:inverter_topologies_controls} summarize inverter model constraints categorized by topology and control laws. From these tables one can derive the feasible set of each inverter topology. 

\begin{table}
    \centering 
    \caption{Inverter Models ($3\times2$) Compared in This Work.}
    \label{table:inverer_constraints}
    \begin{tabular}{l l  l }
        \toprule
         & 4-Leg & 3-Leg \\
         & Section \ref{sec:GFL} & Section \ref{sec:GFL_3w} \\
        \midrule
        \multirow{3}{4.5em}{\textsc{GFL Inverter}}
        & \multirow{3}{11em}{
        \eqref{eq:Ig_sum_0} +\eqref{eq:Ig_limit} + \eqref{eq:Ig_shunt} + \eqref{eq:Ig_Ish_Ii} + \eqref{eq:Si} + \eqref{eq:Sgint} + \eqref{eq:Sgint_Si} + \eqref{eq:Si_pn} +\eqref{eq:Sgint_sum}} 
        & \multirow{3}{11em}{GFL 4-leg model over the phases in $\mathcal{P}$, + \eqref{eq:Ign_0}  }  
        \\ 
        &
        & 
        \\
        \\ \midrule 
        \multirow{3}{4.5em}{\textsc{GFM Inverter}} 
        & \multirow{3}{11em}{\textsc{GFL} 4-Leg model + \eqref{eq:GFM_Umag} + \eqref{eq:GFM_Uangle1}} 
        & \multirow{3}{11em}{\textsc{GFM} 4-Leg model over the phases in $\mathcal{P}$, + \eqref{eq:Ign_0}  }  
        \\
        &
        & 
        \\
        \\ \midrule
        \multirow{2}{4.5em}{\textsc{Solar Inverter}}& \multirow{2}{11em}{\textsc{GFL}/\textsc{GFM} 4-Leg model + \eqref{eq:Pgint_pprimary} + \eqref{eq:Qgint_slack}} &
        \multirow{2}{11em}{\textsc{GFL}/\textsc{GFM} 3-Leg model + \eqref{eq:Pgint_pprimary} + \eqref{eq:Qgint_slack}} 
        \\
        &
        & 
        \\
   \bottomrule
    \end{tabular}
\end{table}

\begin{table}
    \centering 
    \caption{Mapping of Inverter Control Laws to Inverter Models.}
    \label{table:inverter_topologies_controls}
    \def\arraystretch{1.1}
    \begin{tabular}{l l l l l l}
        \toprule
        Topology &  Volt-var & Volt-Watt  & PF & Droop & Unbal. \\
        \midrule
        \textsc{GFL} 4-Leg & 
        \multirow{3}{4.7em}{
        \eqref{eq:control_4wire_Qan},
        \eqref{eq:control_3wire_Qab}, 
        \eqref{eq:control_3wire_Qab_avg}, 
        \eqref{eq:control_3wire_Qag_avg},
        \eqref{eq:control_ODSS_Qag}} & 
        \multirow{3}{4.7em}{
        \eqref{eq:control_4wire_Pan}, 
        \eqref{eq:control_3wire_Pab}, 
        \eqref{eq:control_3wire_Pab_avg}, 
        \eqref{eq:control_3wire_Pag_avg},
        \eqref{eq:control_ODSS_Pag}} & \eqref{eq:control_const_pf} & \quad-- & \multirow{2}{3.1em}{\eqref{eq:control_Ig_012_target} + \eqref{eq:Ig_limit_012}} \\
        & & & & & \\
        & & & & & \\ \midrule
        \textsc{GFL} 3-Leg & 
        \multirow{2}{5.2em}{
        \eqref{eq:control_4wire_Qan},
        \eqref{eq:control_3wire_Qab}, 
        \eqref{eq:control_3wire_Qab_avg}, 
        \eqref{eq:control_3wire_Qag_avg}} & 
        \multirow{2}{5.2em}{
        \eqref{eq:control_4wire_Pan},
        \eqref{eq:control_3wire_Pab}, 
        \eqref{eq:control_3wire_Pab_avg}, 
        \eqref{eq:control_3wire_Pag_avg}} & \eqref{eq:control_const_pf} & \quad-- & \multirow{2}{3.1em}{\eqref{eq:control_Ig_012_target} + \eqref{eq:Ig_limit_012}}  \\
        & & & & & \\ \midrule
        \textsc{GFM} 4-Leg & \quad-- & \quad-- & \ -- & \ \eqref{eq:control_droop} & \ \multirow{2}{3.1em}{\eqref{eq:control_Ig_012_target} + \eqref{eq:Ig_limit_012}} \\
        \\
         \midrule
        \textsc{GFM} 3-Leg & \quad-- & \quad-- & \ -- & \ \eqref{eq:control_droop} & \ \multirow{2}{3.1em}{\eqref{eq:control_Ig_012_target} + \eqref{eq:Ig_limit_012}}    \\
        \\
        \bottomrule
    \end{tabular}
\end{table}



\section{Proof-of-concept studies} 
\label{sec_studies}
We demonstrate how inverter model topologies and controls impact their setpoints and operation in power networks.
We implement the mathematical models in \textsc{JuMP} \cite{lubin_jump_2023} and solve the resulting models using the \textsc{Ipopt} solver \cite{Wachter2006} to a tolerance of 1E-8.
The \textsc{OpenDSS} power flow tolerance is set to 1E-8.

\begin{figure}[tbh]
    \centering
    \includegraphics[width=0.4\columnwidth]{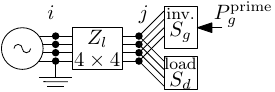}
    \caption{Electrical diagram of 2-bus network to illustrate differences in feasibility between inverter models.}  
    \label{fig:2bus}
\end{figure}

For the upcoming studies, the 3- and 4-leg inverter and unbalanced two-phase load parameters in Table~\ref{table:case_study_data} are used w.r.t. the 2-bus system depicted in Fig. \ref{fig:2bus}.

\begin{table}
    \centering 
    \caption{Case Study Inverter and Load Data.}
    \label{table:case_study_data}
    \def\arraystretch{1.1}
    \begin{tabular}{l l }
        \toprule
        Inverter & Load \\
        \begin{tabular}{l}
             \hline
             ${S}_g^{\text{rating}} = 40 $ kVA  \\
             $I^{\text{rating}}_{g,p} = 52 $ A  \\
            $R_{g,p}=0.015 \ \Omega$ \\
            $X_{g,p}=0.132 \ \Omega$  \\
            $B_{g,p}=1.04$E-7\ S   \\
            $D^{\text{p}}=0.004$ V/W  \\
            $D^{\text{q}}=0.016$ V/var 
        \end{tabular}
         &
        \begin{tabular}{l}
            \hline
            $\mathbf{P}_d = \begin{bmatrix}
                9.0 \\ 4.5 \\ 0.0 
            \end{bmatrix}$ kW
            \\
            \\
            $\mathbf{Q}_d = \begin{bmatrix}
                4.36 \\ 2.18 \\ 0.0
            \end{bmatrix}$ kvar
        \end{tabular}
            \\
        \bottomrule
    \end{tabular}
\end{table}

\subsection{Validation w.r.t. \textsc{OpenDSS}}
We first validate our implementation w.r.t. \textsc{OpenDSS}. GFL inverters of 3-, 4-leg and 1-leg designs with Volt-var control modes \eqref{eq:control_3wire_Qag_avg}, \eqref{eq:control_3wire_Qan_avg}, and \eqref{eq:control_ODSS_Qag}, shown in Table~\ref{tab:vv_vw}, are modelled and the comparison reported in Table~\ref{table:validation_GFL_3w_4w} shows the accuracy of the models with respect to voltage magnitudes and reactive power generation. Note that to solve the control laws in \textsc{OpenDSS}, we tighten var and voltage convergence tolerances (\texttt{VarChangeTolerance}, \texttt{ VoltageChangeTolerance}), compared to the default values, and set them to 1E-4 and 1E-5, respectively, which improves the accuracy of the control setpoints to some extent. However, this reveals that these tolerances are the most likely reasons for reduced accuracy, which is the source of the mismatch in ${Q}_{g}$. 

In order to validate the results for the averaged phase-to-neutral Volt-var control \eqref{eq:control_3wire_Qan_avg}, the \textsc{OpenDSS} case file had to be adapted by adding a lossless transformer between the bus $j$ and the load/inverter in Figure~\ref{fig:2bus}. The neutral at the secondary of the transformer is then grounded and the voltage and reactive power at the transformer primary bus, bus $j$, is compared to our optimization model results.

To validate the 4-leg GFL inverters with control function \eqref{eq:control_ODSS_Qag}, we replaced the 4-leg inverter in Figure~\ref{fig:2bus} by three single-phase inverters in OpenDSS connected to each phase with ${S}_g^{\text{rating}} = 40/3 $ kVA and $P_g = 35/3$ kW. Due to \textsc{OpenDSS} being very sensitive to var and voltage convergence tolerances in this case, these values were set higher to 1E-2 and 1E-3 respectively, which is the main cause of the inaccurate voltage and reactive power generation w.r.t our implementation. Note that \textsc{OpenDSS} inaccuracy is confirmed by validating if $(\mathbf{U}^{\text{mag}}_{g,a}, \mathbf{Q}_{g})$ is on the $f^{VV}$ curve.

\begin{table}
    \centering 
    \caption{Validating inverter model against OpenDSS: GFL inverter with Volt-var controls \eqref{eq:control_3wire_Qag_avg}, \eqref{eq:control_3wire_Qan_avg}, and \eqref{eq:control_ODSS_Qag}}
    \label{table:validation_GFL_3w_4w}
    \def\arraystretch{1.1}
    \begin{tabular}{l  l l }
        \toprule
        &  \textsc{OpenDSS} & Our implementation \\
        \midrule
        \multirow{2}{5em}{\textsc{GFL} 3-Leg  $f^{VV}$ \eqref{eq:control_3wire_Qag_avg}} & 
        $\mathbf{Q}_{g} =\left[\begin{smallmatrix}-0.3965 \\ -0.3965 \\ -0.3965 \end{smallmatrix}\right]$  kvar& 
        $\mathbf{Q}_{g} = \left[\begin{smallmatrix} -0.3962 \\ -0.3962 \\ -0.3962  \end{smallmatrix}\right]$ kvar\\
        & 
        $\mathbf{U}^{\text{mag}}_{g} = \left[\begin{smallmatrix} 
            0.9765 \\
            1.0048 \\
            1.0278 \\
            0.0736
        \end{smallmatrix}\right]$ pu & 
        $\mathbf{U}^{\text{mag}}_{g} = \left[\begin{smallmatrix} 
            0.9765 \\
            1.0048 \\
            1.0278 \\
            0.0736
        \end{smallmatrix}\right]$ pu
        \\ \midrule
        \multirow{2}{5em}{\textsc{GFL} 4-Leg   $f^{VV}$ \eqref{eq:control_3wire_Qag_avg}} & 
        $\mathbf{Q}_{g}= \left[\begin{smallmatrix}-0.3946 \\ -0.3946 \\ -0.3946 \end{smallmatrix}\right]$ kvar & 
        $\mathbf{Q}_{g}= \left[\begin{smallmatrix}-0.3943 \\ -0.3943 \\ -0.3943 \end{smallmatrix}\right]$ kvar\\
        & 
        $\mathbf{U}^{\text{mag}}_{g} = \left[\begin{smallmatrix} 
            0.9748 \\
            1.0123 \\
            1.0220 \\
            0.0615
        \end{smallmatrix}\right]$ pu & 
        $\mathbf{U}^{\text{mag}}_{g} = \left[\begin{smallmatrix} 
            0.9748 \\
            1.0123 \\
            1.0220 \\
            0.0615
        \end{smallmatrix}\right]$ pu
        \\ \midrule
        \multirow{2}{5em}{
\textsc{GFL} 4-Leg   $f^{VV}$ \eqref{eq:control_3wire_Qan_avg}} & 
        $\mathbf{Q}_{g}= \left[\begin{smallmatrix}-0.4872 \\ -0.4872\\ -0.4872 \end{smallmatrix}\right]$ kvar & 
        $\mathbf{Q}_{g}= \left[\begin{smallmatrix}-0.4877 \\ -0.4877 \\ -0.4877 \end{smallmatrix}\right]$ kvar \\
        & 
        $\mathbf{U}^{\text{mag}}_{g} = \left[\begin{smallmatrix} 
             0.9741 \\
            1.0117 \\
            1.0213 \\
            0.0617
        \end{smallmatrix}\right]$ pu & 
        $\mathbf{U}^{\text{mag}}_{g} = \left[\begin{smallmatrix} 
            0.9741 \\
            1.0117 \\
            1.0213 \\
            0.0617
        \end{smallmatrix}\right]$ pu
        \\ \midrule
        \multirow{2}{5em}{\textsc{GFL} 4-Leg   $f^{VV}$ \eqref{eq:control_ODSS_Qag}} & 
        $\mathbf{Q}_{g}= \left[\begin{smallmatrix}\hfill 1.8737 \\ -1.8959 \\ -1.8849 \end{smallmatrix}\right]$ kvar & 
        $\mathbf{Q}_{g}= \left[\begin{smallmatrix}\hfill 1.8986 \\ -1.8657 \\ -1.8582 \end{smallmatrix}\right]$ kvar \\
        & 
        $\mathbf{U}^{\text{mag}}_{g} = \left[\begin{smallmatrix} 
             0.9851 \\
             1.0142 \\
             1.0142 \\
             0.0
        \end{smallmatrix}\right]$ pu & 
       $\mathbf{U}^{\text{mag}}_{g} = \left[\begin{smallmatrix} 
            0.9853 \\
            1.0144 \\
            1.0144 \\
            0.0
        \end{smallmatrix}\right]$ pu
            \\
        \bottomrule
    \end{tabular}
\end{table}

\subsection{Base Case: GFL 4-Leg vs. 3-Leg Set Points}
The inverter topology significantly impacts feasible operating points even for a basic objective such as network loss minimization. Table~\ref{table:GFL_3w_4w_minloss} shows how GFL power and current values are different for 3 vs. 4-leg  while the voltage magnitude at the PCC is close for both cases. With this observation, in the next sections we investigate behaviour of each topology for GFL and GFM inverters and for different control objectives.
To minimize the network power losses, we use, 
\begin{IEEEeqnarray}{C}
    \label{eq_min_loss}
      \min \sum_l P^{\text{loss}}_{l} \text{ with } P^{\text{loss}}_{l} = \sum_{p \in \mathcal{P}} P_{lij,p} + P_{lji,p}.
\end{IEEEeqnarray}
Minimization of negative sequence current flowing through lines is, 
\begin{IEEEeqnarray}{C}
    \label{eq:IUF2}
      \min \sum_{lij} |I_{lij,2}|.
\end{IEEEeqnarray}
Minimization of generation cost is,
\begin{IEEEeqnarray}{C}
    \label{eq:gencost}
      \min \sum_{g} c_g  \sum_{p \in \mathcal{P}} \mathbf{P}_{g}.
\end{IEEEeqnarray}
We set up a simple case study with one branch, depicted in Fig. \ref{fig:2bus}.

\begin{table}
    \centering 
    \caption{3-Leg and 4-Leg GFL inverter set point comparison}
    \label{table:GFL_3w_4w_minloss}
    \begin{tabular}{l c c c}
        \toprule
         & $\mathbf{U}_{g} $ (pu) & $\Sgint $ (kW+$j$kvar) & $\mathbf{I}_g$ (A) \\
        \midrule
        3-Leg & 
        $\left[\begin{smallmatrix} 
            0.9270 \angle -2.33 \hfill \\
            0.9912 \angle -121.56 \hfill \\
            0.9948 \angle 119.91 \hfill \\
            0.0880 \angle 14.82 \hfill
        \end{smallmatrix}\right]$ & 
        $\left[\begin{smallmatrix} 
            1.988 + \imagnumber 0.302 \\
            0.0  \hfill+ \imagnumber 0.901 \\
            0.618 - \imagnumber 1.2 \hfill \\
        \end{smallmatrix}\right]$ &
        $\left[\begin{smallmatrix} 
             \hfill 2.13  \hfill - \imagnumber 0.413 \\
            -0.775 + \imagnumber 0.476 \\
            -1.355 - \imagnumber 0.063 \\
        \end{smallmatrix}\right]$
        \vspace{1ex} 
        \\ \midrule
        4-Leg & $\left[\begin{smallmatrix} 
            0.9295 \angle -2.44 \hfill \\
            0.9889 \angle -121.53 \hfill \\
            0.9937 \angle 119.72 \hfill \\
            0.0813 \angle 16.01 \hfill
        \end{smallmatrix}\right]$ & 
        $\left[\begin{smallmatrix} 
            1.675 + \imagnumber 0.415 \\
            0.0  \hfill+ \imagnumber 0.764 \\
            0.320 - \imagnumber 1.195 \\
            0.0  \hfill+ \imagnumber 0.0 \hfill
        \end{smallmatrix}\right]$ & 
        $\left[\begin{smallmatrix} 
            \hfill 1.924 - \imagnumber 0.628 \\
            -0.599 + \imagnumber 0.412 \\
            -1.150 - \imagnumber 0.4 \hfill \\
            -0.175 + \imagnumber 0.616
        \end{smallmatrix}\right]$
        \\
        \bottomrule
    \end{tabular}
\end{table}

\subsection{Phase Unbalance Compensation: {GFL} 4-Leg vs. 3-Leg}
We demonstrate phase unbalance compensation capability of each inverter topology with the objective of minimizing the grid current negative sequence component \eqref{eq:IUF2}.
We consider two scenarios where the inverter's primary source power is 1) higher and 2) lower than the total load. Results are shown in Table~\ref{table:unbalance_compensation_GFL_3w_4w}.

The load conductor currents, shown in the right column, show the phase unbalance of the two-phase load, with the neutral conductor current being an indication of the zero-sequence component. 
Note that in all scenarios in Table~\ref{table:unbalance_compensation_GFL_3w_4w}, the load active and reactive powers are the same as in Table~\ref{table:case_study_data}, however, the load currents that are shown in the figures are not necessarily the same, hence different magnitudes of current phasors with respect to the circle.

For the first scenario shown in the top two rows, while both inverters are capable of delivering negative-sequence current, the 3-leg inverter is not capable of delivering zero-sequence current which in turn must drawn from the grid. 
The 4-leg inverter, however, can deliver both zero- and negative-sequence currents, and the reason for small zero-sequence current from the grid is the objective function which minimizes only the negative-sequence component. These results are consistent with the \emph{simulations} of Stuyts et al. \cite{Stuyts2016}, but were now derived using a constrained nonlinear \emph{optimization} algorithm.

In the next scenario we test inverter behaviour when the available primary source power is too small to supply the load by itself. 
In this case, the 4-leg inverter is still able to minimize the negative sequence, while the 3-leg inverter is limited by how much it can compensate for phase currents through only the phase conductors.

\begin{table}
    \centering 
    \caption{3-Leg and 4-Leg GFL inverter unbalance compensation comparison for negative sequence current minimization \eqref{eq:IUF2}. Color guide: $I_a$: blue, $I_b$: red, $I_c$: green, $I_n$: black.}
    \label{table:unbalance_compensation_GFL_3w_4w}
    \begin{tabular}{l  l }
        \toprule
        Topology & \ \ \ \ \ \ Grid  \ \ \  \ \ \ + \ \ \ \ \ \
        Inverter \ \ \ \ \ = \ \ \ \ \ \ 
        Load \\
        \midrule
        \\
        \begin{tabular}{l}
            \textsc{GFL} 4-Leg
        \end{tabular} & 
        \hspace{-7ex}
        \begin{tabular}{l}
        \includegraphics[width=0.75\columnwidth]{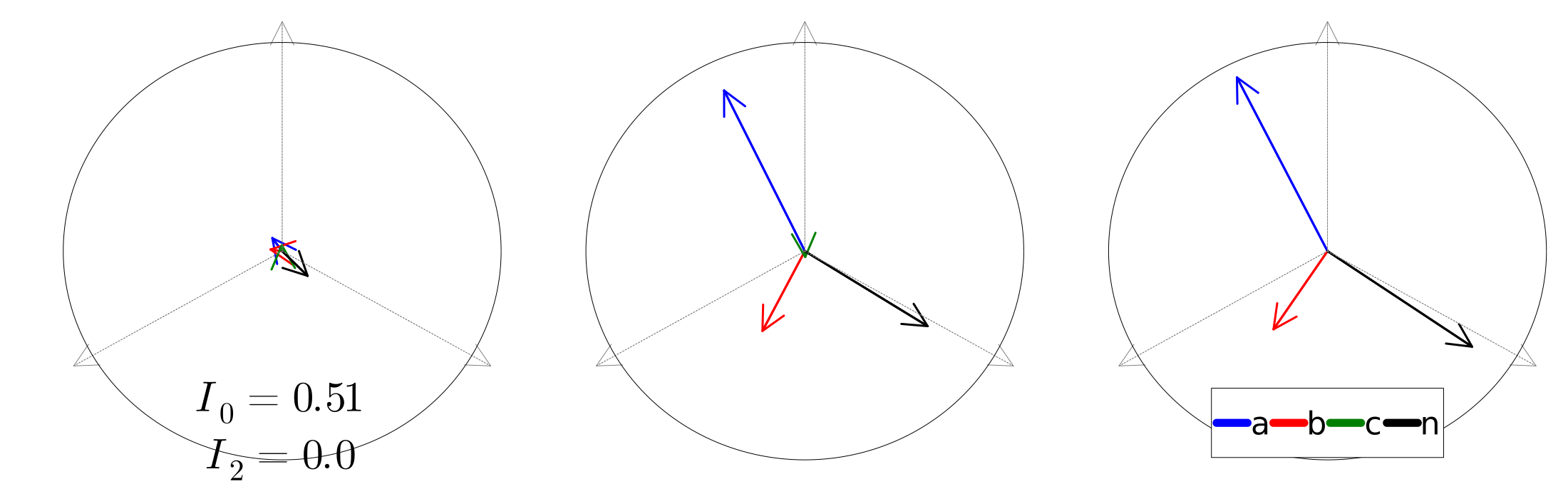}
        \end{tabular}
        \\
        \begin{tabular}{l}
            \textsc{GFL} 3-Leg
        \end{tabular} & 
         \hspace{-7ex}
        \begin{tabular}{l}
        \includegraphics[width=0.75\columnwidth]{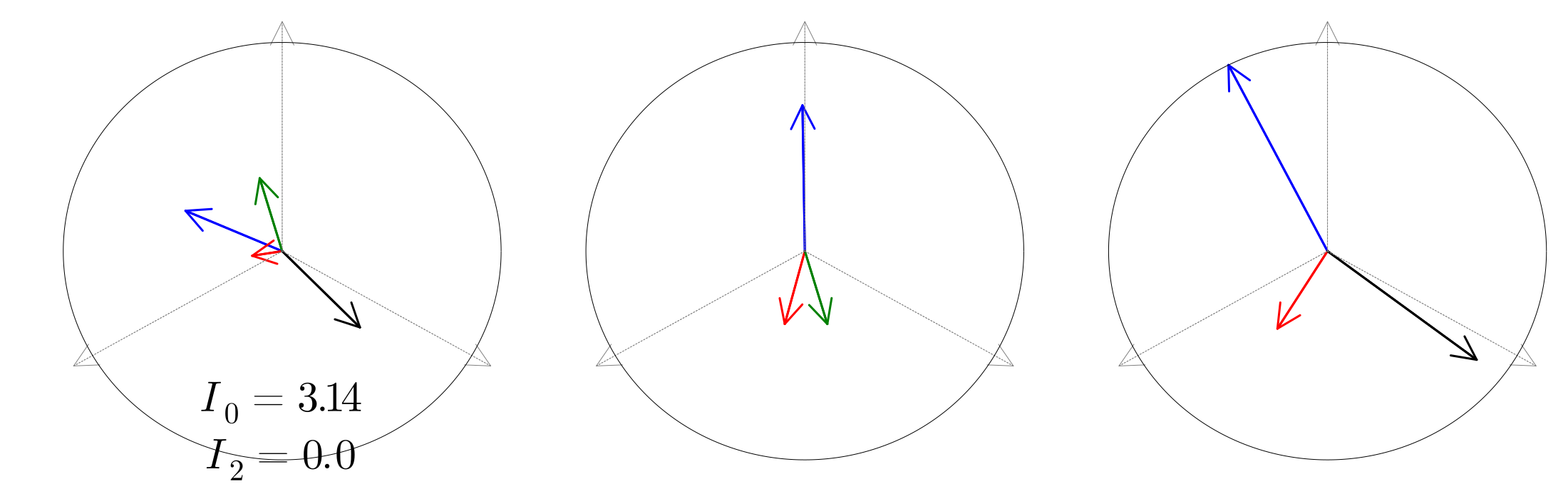}
        \end{tabular}
        \\
        \hline
        \begin{tabular}{l}
            \textsc{GFL} 4-Leg \\ ${P}^{\text{prime}}_{g} \le 3$
        \end{tabular} & 
         \hspace{-7ex}
        \begin{tabular}{l}
        \includegraphics[width=0.75\columnwidth]{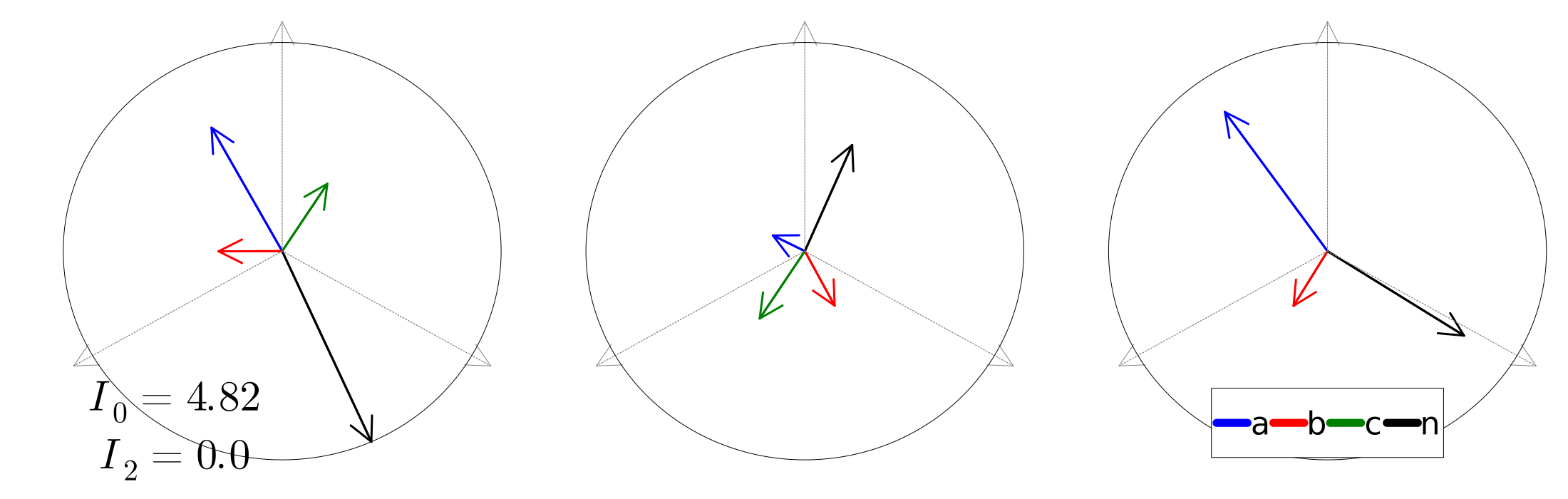}
        \end{tabular}
        \\
        \begin{tabular}{l}
            \textsc{GFL} 3-Leg \\ ${P}^{\text{prime}}_{g} \le 3$
        \end{tabular} & 
         \hspace{-7ex}
        \begin{tabular}{l}
        \includegraphics[width=0.75\columnwidth]{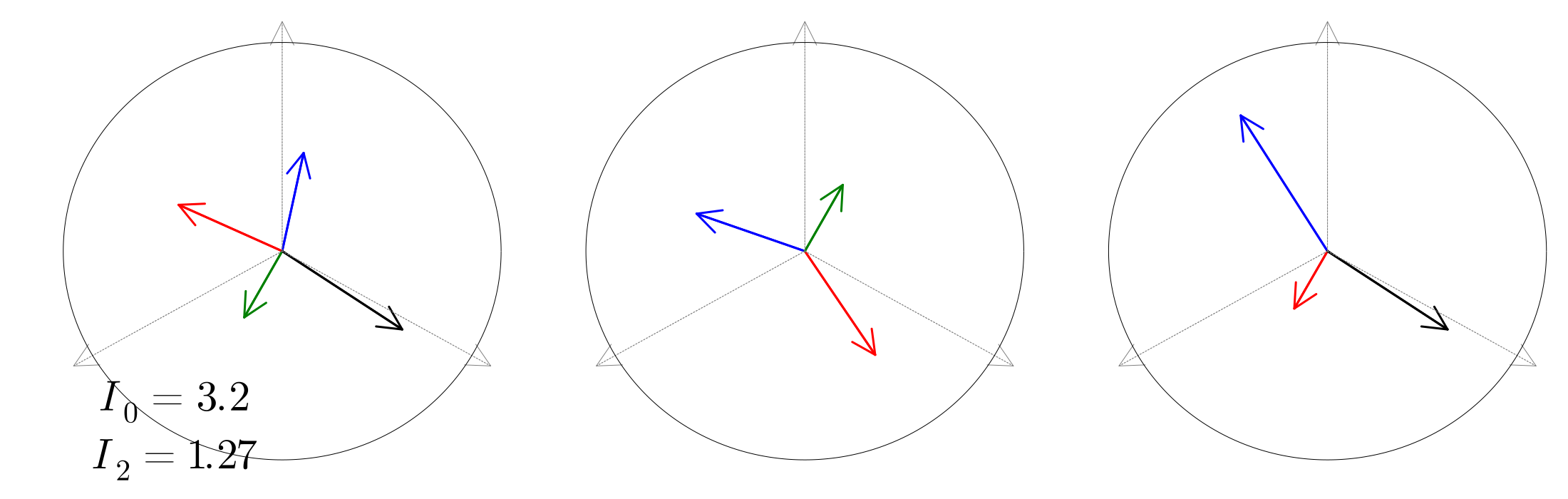}
        \end{tabular}
        \\
        \bottomrule
    \end{tabular}
\end{table}

\subsection{GFL Volt-var and Volt-Watt Control Laws}
\begin{figure}[tbh]
    \centering
    \includegraphics[width=1\columnwidth]{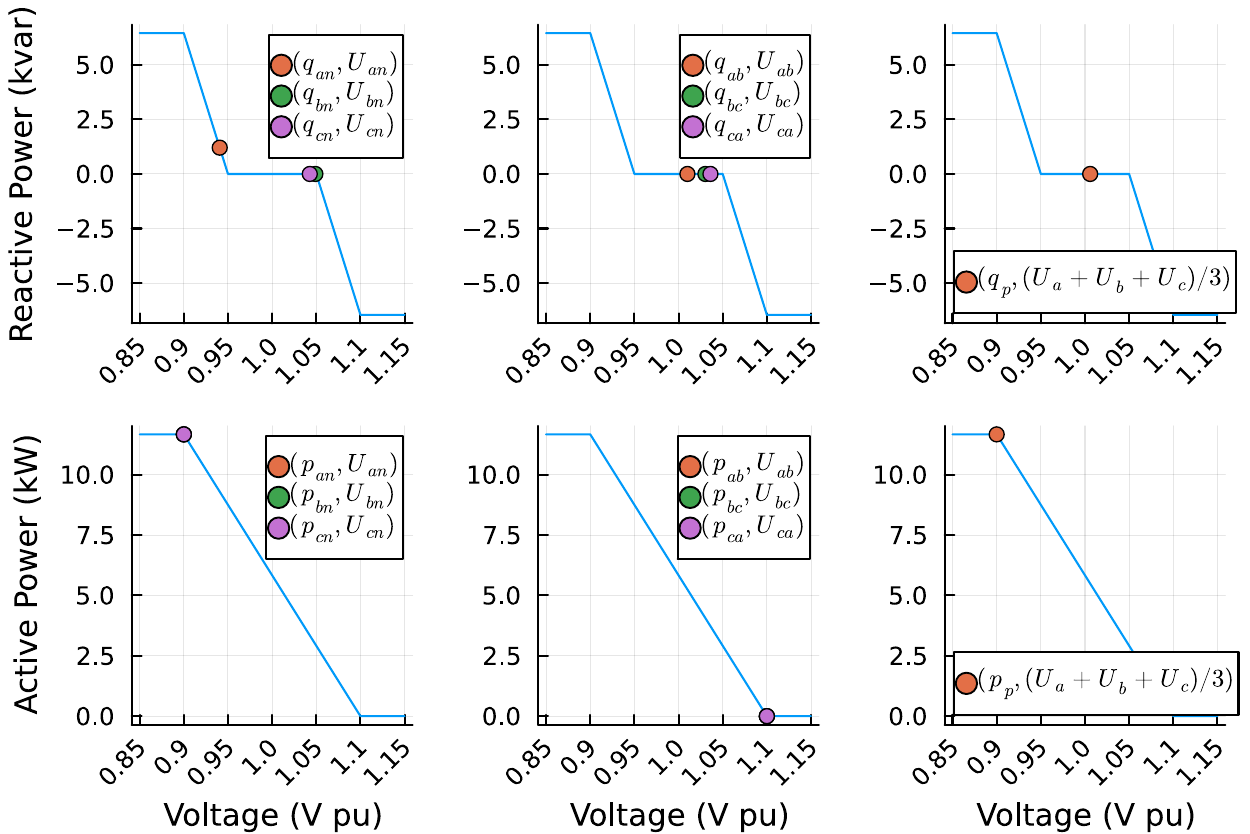}
    \caption{Volt-var and Volt-Watt control laws in 4-Leg GFL inverters. From left to right 1) phase-to-neutral control, 2) phase-to-phase control and 3) averaged phase-to-ground control.}  
    \label{fig:vv_vw_4w_GFL}
\end{figure}

We next demonstrate Volt-var and Volt-Watt control laws in GFL inverters, based on Australian gridcode standards \cite{powerquality_ENA}, for phase-to-neutral, and phase-to-phase, and averaged phase-to-ground voltage magnitudes. Fig.~\ref{fig:vv_vw_4w_GFL} shows how each control law impacts voltage for the same case study with a 4-leg GFL inverter. It is observed that the phase-to-neutral and phase-to-phase control laws support unbalanced phase power flows, whereas the control laws as a function averaged phase-to-ground voltages balance the power delivery through the phase conductors.

\subsection{GFM and {GFL} 4-Leg Inverters and Relaxation Hierarchy}

GFM inverters enforce the voltage phase angle differences at $2\pi/3$, but the  voltage magnitudes can 1) be constrained to a \emph{set point}, e.g. 1.0 pu, 2) be only \emph{equal} for all the phases as in \eqref{eq:GFM_Umag}, or 3) be determined from a \emph{droop} control law \eqref{eq:control_droop}. 
Table~\ref{table:GFM_GFL_4w_relaxations} shows the internal voltage vector for each GFM control mode where the voltage angle differences are $2\pi/3$, but the voltage magnitudes vary according to the control mode.

Furthermore, the GFL inverter is considered as a relaxation of the GFM inverter with the internal voltage magnitudes and angles are linked.
Table~\ref{table:GFM_GFL_4w_relaxations} provides cost, grid current negative-sequence component, and inverter internal voltage in the top-down hierarchical order where the GFM 4-w with droop control is the most constrained form and the GFL 4-w is the most relaxed form of the 4-leg inverters, the ranking that is also reflected by the top-to-bottom cost reduction the cost minimization problem \eqref{eq:gencost}.

In the negative-sequence current compensation case, the GFL 4-leg is able to deliver more negative-sequence current compared to the equal-voltage magnitudes GFM. The droop-controlled GFM has the lowest current negative-sequence compensation capability. 

\begin{table}
    \centering 
    \caption{Cost minimization for GFM and \textsc{GFL} 4-Leg model.}
    \label{table:GFM_GFL_4w_relaxations}
    \def\arraystretch{1.1}
    \begin{tabular}{l c c c}
        \toprule
         & Cost ($\$/h$)  & $|I_{lij,2}|$ & $\Ugint$   \\
        \midrule
        \textsc{GFM} 4-Leg droop \eqref{eq:control_droop} &
        136.17 & 1.78 &
        $\left[\begin{smallmatrix}
            0.968 \angle 0.3  \hfill \\
            0.997 \angle -119.7  \hfill \\
            1.014 \angle 120.3  \hfill \\
            0.0 \angle 0  \hfill
        \end{smallmatrix}\right]$ 
        \\ \midrule
        \textsc{GFM} 4-Leg set point &
        135.69 & 0.79 &
        $\left[\begin{smallmatrix}
            1.0 \angle 0.3 \hfill     \\
            1.0 \angle -119.7  \hfill \\
            1.0 \angle 120.3  \hfill \\
            0.0 \angle 0  \hfill
        \end{smallmatrix}\right]$
        \\ \midrule
        \textsc{GFM} 4-Leg equal \eqref{eq:GFM_Umag} &
        135.66 & 0.76 &
         $\left[\begin{smallmatrix}
            0.955 \angle 0.3  \hfill \\
            0.955 \angle -119.7  \hfill \\
            0.955 \angle 120.3  \hfill \\
            0.0 \angle 0  \hfill
        \end{smallmatrix}\right]$
        \\ \midrule
        \textsc{GFL} 4-Leg &
        135.52 & 0.22 & -
        \\
        \bottomrule
        \\
    \end{tabular}
\end{table}

\subsection{Optimal Power Flow with GFM and GFL Inverter Models}
We showcase operation of GFM and GFL inverters with 3 and 4-leg topologies in a 4-wire (i.e. $4\times 4$ line impedance matrices) adaptation of the European low voltage test Feeder \cite{EPRI8063903} with 907 buses. 
We consider three scenarios where we have 6 DERs connected to the grid a) without inverters, b) with GFL 4-leg inverters, and c) with 3 GFM 4-leg, 2 GFL 4-leg and 1 GFL 3-leg inverters. 
The objective is cost minimization in the network and the voltage negative-sequence component is limited to $2 \%$ as in \eqref{eq:Ui2_limit}.

\begin{figure}[tbh]
    \centering
    \includegraphics[width=1\columnwidth]{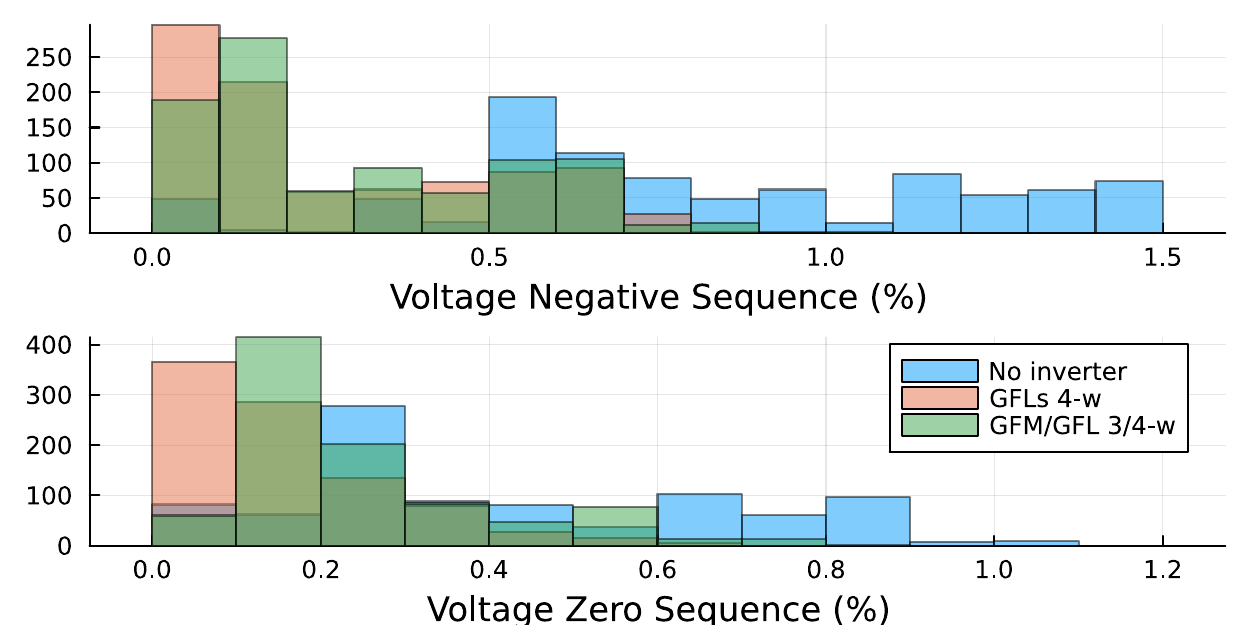}
    \caption{Negative and zero sequence components of network voltages.}  
    \label{fig:opf_v0_v2}
\end{figure}

The results Table~\ref {tab:OPF_GFL_GFM} show objective value as well as solve time and \textsc{Ipopt} iterations. The objective values show that the most relaxed scenario without inverter models has the lowest cost, whereas the most constrained scenario with a mix of GFM and GFL inverters has the highest cost. The solve time and iterations also support this relaxation hierarchy as both increase with  more constraints added to the network. We also show the worst voltage negative-sequence in the grid. 
Fig.~\ref{fig:opf_v0_v2} shows voltage zero- and negative-sequence components throughout the network for the three scenarios. 

\begin{table}[b]
\caption{Optimal power flow results for EU LV test feeder with a voltage source and 6 inverter-based dispatchable resources. }
\label{tab:OPF_GFL_GFM}
\centering
\def\arraystretch{1.1}
\begin{tabular}{l l l l l } 
    \toprule    
    Inverter model & \multirow{2}{3.5em}{Objective (\$/h)} & \multirow{2}{3.8em}{Solve time (s)} & \multirow{2}{3.5em}{\textsc{Ipopt} iterations} & \multirow{2}{5em}{Worst $|U_{i,2}|$ (\%)}
    \\
    &  & 	&  &  \\
    \midrule
    Original: no Inverter & 69.999  & 2.59 & 22  &  1.65 \\ 
    \textsc{GFL} 4-Leg & 69.656  &  2.30 &  20 & 0.95  \\ 
    \textsc{GFM}/\textsc{GFL} 3/4-Leg & 69.934  & 2.47	& 20  &  1.01 \\
    \bottomrule
\end{tabular}
\end{table}


\section{Conclusions} 
\label{sec_conclusions}
A mathematical framework to represent 3 and 4-leg inverters with various levels of control functionalities was developed.
We limit ourselves (nonlinear) algebraic to enable the use of in numerical optimization contexts.
Significant differences in feasibility and dispatch costs between the different inverter models were illustrated numerically, and justified from the perspective of relaxation/restriction.
Finally, a case study integrating these models in a four-wire OPF engine \cite{CLAEYS2022108522} is presented, and we did not observe significant changes in convergence despite the increased level of detail in the model. 
The grid-forming inverter models are likely particularly useful in the context of microgrid-focused studies. 
The grid-following ones enable more detailed studies on battery dispatch optimization in LV networks with PV with distributed Volt-var/Watt control. 

We observe that most published studies on inverter dispatch optimization can be understood as representing only 4-leg grid-following inverters, with or without a power sharing constraint. 3-leg inverter models are restrictions of such models, so real-world costs are expected to be underestimated.

Future work include convex relaxation models for inverters, and studying the implications of control degrees of freedom on the sizing of battery storage systems or active filters, and extensions in the direction of harmonic (optimal) power flow. We also extend the inverter models to include various prime mover models, e.g. stationary or EV batteries, and solar PV. 
Establishing inverter models from data-driven characterization through hardware-in-the-loop testing is a direction for future work.

\bibliographystyle{IEEEtran}
\bibliography{PSCCfullpaper}


\end{document}

%% file: newcommands.tex
\newcommand{\indexnode}{i}
\newcommand{\indexconductor}{p}
\newcommand{\indextime}{t}
\newcommand{\indexconverter}{c}
\newcommand{\indexstorage}{s}

\newcommand{\setconductor}{\mathcal{P}}
\newcommand{\setTime}{\mathcal{T}}

\newcommand{\symtime}{T}
\newcommand{\symenergy}{E}
\newcommand{\symcharge}{Q}

\newcommand{\symcomplexpower}{S}
\newcommand{\symactivepower}{P}
\newcommand{\symreactivepower}{Q}
\newcommand{\symcurrent}{I}
\newcommand{\symcurrentlifted}{L}
\newcommand{\symvoltage}{U}
\newcommand{\symvoltagelifted}{W}
\newcommand{\symimpedance}{Z}
\newcommand{\symresistance}{R}
\newcommand{\symreactance}{X}
\newcommand{\symstatus}{s}
\newcommand{\symbinary}{z}
\newcommand{\symefficiency}{\eta}

\newcommand{\paramcolor}{black}
\newcommand{\imagnumber}[0]{\textcolor{\paramcolor}{j}}

\newcommand{\imag}{j}

\newcommand{\ab}{\alpha \beta \gamma}

\newcommand{\Ugint}{\mathbf{U}^{\text{int}}_{g}}
\newcommand{\Ugintpn}{\mathbf{U}^{\text{pn,int}}_{g}}
\newcommand{\Ugintab}{\mathbf{U}^{\text{int}}_{g, \alpha \beta \gamma}}
\newcommand{\Ugintpnab}{\mathbf{U}^{\text{pn,int}}_{g, \alpha \beta \gamma}}
\newcommand{\Uginta}{\mathbf{U}^{\text{int}}_{g, \alpha}}
\newcommand{\Ugintb}{\mathbf{U}^{\text{int}}_{g, \beta}}

\newcommand{\Sgint}{\mathbf{S}^{\text{int}}_{g}}
\newcommand{\Pgint}{\mathbf{P}^{\text{int}}_{g}}
\newcommand{\Qgint}{\mathbf{Q}^{\text{int}}_{g}}

\newcommand{\Sgintpn}{\mathbf{S}^{\text{pn,int}}_{g}}
\newcommand{\Sin}{\mathbf{S}^{\text{pn}}_{i}}

\newcommand{\Scrating}{\textcolor{\paramcolor}{\symcomplexpower_{\indexconverter}^{\text{rating}}}}
\newcommand{\Icrating}{\textcolor{\paramcolor}{\symcurrent_{\indexconverter}^{\text{rating}}}}

\newcommand{\Scratingp}{\textcolor{\paramcolor}{\symcomplexpower_{\indexconverter, \indexconductor}^{\text{rating}}}}
\newcommand{\Icratingp}{\textcolor{\paramcolor}{\symcurrent_{\indexconverter, \indexconductor}^{\text{rating}}}}

\newcommand{\statusct}{\textcolor{\paramcolor}{\symstatus_{\indexconverter, \indextime}}}

\newcommand{\Sct}{\symcomplexpower_{\indexconverter, \indextime}}
\newcommand{\Pct}{\symactivepower_{\indexconverter, \indextime}}
\newcommand{\Qct}{\symreactivepower_{\indexconverter, \indextime}}

\newcommand{\Scpt}{\symcomplexpower_{\indexconverter, \indexconductor, \indextime}}
\newcommand{\Pcpt}{\symactivepower_{\indexconverter, \indexconductor, \indextime}}
\newcommand{\Qcpt}{\symreactivepower_{\indexconverter, \indexconductor, \indextime}}

\newcommand{\Pcat}{\symactivepower_{\indexconverter, a, \indextime}}
\newcommand{\Pcbt}{\symactivepower_{\indexconverter, b, \indextime}}
\newcommand{\Pcct}{\symactivepower_{\indexconverter, c, \indextime}}

\newcommand{\Qcat}{\symreactivepower_{\indexconverter, a, \indextime}}
\newcommand{\Qcbt}{\symreactivepower_{\indexconverter, b, \indextime}}
\newcommand{\Qcct}{\symreactivepower_{\indexconverter, c, \indextime}}

\newcommand{\Ict}{\symcurrent_{\indexconverter, \indextime}}
\newcommand{\Lct}{\symcurrentlifted_{\indexconverter, \indextime}}

\newcommand{\Icpt}{\symcurrent_{\indexconverter,  \indexconductor,\indextime}}
\newcommand{\Lcpt}{\symcurrentlifted_{\indexconverter, \indexconductor, \indextime}}

\newcommand{\Uit}{\symvoltage_{\indexnode, \indextime}}
\newcommand{\Zc}{\textcolor{\paramcolor}{\symimpedance_{\indexconverter}}}

\newcommand{\Wit}{\symvoltagelifted_{\indexnode, \indextime}}

\newcommand{\Uipt}{\symvoltage_{\indexnode, \indexconductor, \indextime}}
\newcommand{\Zcp}{\textcolor{\paramcolor}{\symimpedance_{\indexconverter, \indexconductor}}}
\newcommand{\Rcp}{\textcolor{\paramcolor}{\symresistance_{\indexconverter, \indexconductor}}}
\newcommand{\Xcp}{\textcolor{\paramcolor}{\symreactance_{\indexconverter, \indexconductor}}}

\newcommand{\Uiptmin}{\textcolor{\paramcolor}{\symvoltage_{\indexnode, \indexconductor}^{\text{min}}}}
\newcommand{\Uiptmax}{\textcolor{\paramcolor}{\symvoltage_{\indexnode, \indexconductor}^{\text{max}}}}

\newcommand{\Wipt}{\symvoltagelifted_{\indexnode, \indexconductor, \indextime}}

\newcommand{\Pstor}{\symactivepower_{\indexconverter, \indextime}^{\text{stor}}}
\newcommand{\Sstor}{\symcomplexpower_{\indexconverter, \indextime}^{\text{stor}}}

\newcommand{\Pext}{\textcolor{\paramcolor}{\symactivepower_{\indexconverter, \indextime}^{\text{ext}}}}
\newcommand{\Sext}{\textcolor{\paramcolor}{\symcomplexpower_{\indexconverter, \indextime}^{\text{ext}}}}

\newcommand{\Qint}{\textcolor{\paramcolor}{\symreactivepower_{\indexconverter, \indextime}^{\text{int}}}}

\newcommand{\Pcstor}{\symactivepower_{\indexconverter, \indextime}^{\text{c}}}
\newcommand{\Pdstor}{\symactivepower_{\indexconverter, \indextime}^{\text{d}}}

\newcommand{\bc}{\symbinary_{\indexconverter, \indextime}^{\text{c}}}

\newcommand{\Erating}{\textcolor{\paramcolor}{\symenergy_{\indexconverter}^{\text{max}}}}
\newcommand{\Pcstorrating}{\textcolor{\paramcolor}{\symactivepower_{\indexconverter}^{\text{c,max}}}}
\newcommand{\Pdstorrating}{\textcolor{\paramcolor}{\symactivepower_{\indexconverter}^{\text{d,max}}}}

\newcommand{\etacstor}{\textcolor{\paramcolor}{\symefficiency_{\indexconverter}^{\text{c}}}}
\newcommand{\etadstor}{\textcolor{\paramcolor}{\symefficiency_{\indexconverter}^{\text{d}}}}

\newcommand{\Ect}{\symenergy_{\indexconverter, \indextime}}
\newcommand{\Ectinit}{\textcolor{\paramcolor}{\symenergy_{\indexconverter}^{\text{init}}}}
\newcommand{\Ectprev}{\symenergy_{\indexconverter, \indextime-1}}

\newcommand{\Ectone}{\symenergy_{\indexconverter, \indextime=1}}

\newcommand{\Ectend}{\symenergy_{\indexconverter, \indextime=n}}

\newcommand{\Tt}{\textcolor{\paramcolor}{\symtime_{\indextime}}}